\documentclass[12pt,draft]{article}

\usepackage[cp1251]{inputenc}
\usepackage[T2A]{fontenc}
\usepackage[english,russian,ukrainian]{babel}
\usepackage{amsmath,amsfonts,amssymb}
\usepackage{geometry}
\usepackage{cite}
\begin{document}
\begin{flushleft}

\large\textbf{О. В. Дяченко (O. V. Diachenko)}

\medskip

\small{(Національний технічний університет України "'Київський політехнічний інститут імені Ігоря Сікорського"')}

\medskip

\large\textbf{В. М. Лось (V. M. Los)}

\medskip

\small{(Національний технічний університет України "'Київський політехнічний інститут імені Ігоря Сікорського"')}

\large\textbf{Про класичність розв'язків крайової задачі для параболічної системи другого порядку}

\medskip

\large\textbf{On the Classicality of Solutions to a Boundary Value Problem for a Second-Order Parabolic System}

\end{flushleft}

\normalsize

Досліджується параболічна початково-крайова задача для системи двох диференціальних рівнянь з двома крайовими умовами різного порядку~--- Діріхле і Неймана. Вона зустрічається, зокрема, у теорії тепломасообміну. Отримано достатні умови класичності узагальненого розв'язку задачі. Їх сформульовано у термінах належності даних задачі  узагальненим анізотропним просторам Соболєва.

\smallskip

We study a parabolic initial-boundary-value problem for a system of two differential equations with two boundary conditions of different orders,
the Dirichlet and Neumann ones. It occurs specifically in the heat-mass transfer theory. We find sufficient conditions for the generalized solution of the problem to be classical. They are formulated in terms of the belonging of the problem data to generalized anisotropic Sobolev spaces.

\bigskip

\textbf{1. Вступ.}
При дослідженні параболічних задач важливими є питання регулярності, зокрема, класичності їх узагальнених розв'язків.
Під класичним розуміють неперервно диференційовний розв'язок, який  задовольняє задачу в термінах класичних похідних.
Відповідь на ці питання дають, як правило, шляхом формулювання умов приналежності правих частин задачі певним функціональним просторам.
Чим тонше градуйована вибрана шкала функціональних просторів, тим точніший результат буде отримано.
До певного часу параболічні задачі вивчались переважно у функціональних анізотропних просторах Соболєва та Гельдера \cite{Solonnikov65,Ivasyshen90, Eidelman94,ZhitarashuEidelman98,Ilin60, IlinKalashOleinik62}, які параметризуються числами.
В останні десятиріччя активно розвивається дослідження параболічних рівнянь і систем в різних інших шкалах функціональних просторів \cite{DenkHieberPruess07,Lindemulder20,DongKim15,Hummel21,LeCronePrussWilke14}.
Одними з них є шкали узагальнених анізотропних просторів Соболєва \cite{LosMikhailetsMurach21Monograph,LosMikhailetsMurach21CPAA,LosMikhailetsMurach17CPAA}.
Вони параметризуються крім числових показників регулярності ще
додатковим функціональним параметром. Тому їх використання дозволяє отримати більш точні результати, ніж це можливо в межах класичних
шкал Соболєва і Гельдера.
Шкали ізотропних узагальнених просторів Соболєва знайшли численні застосування в теорії еліптичних задач \cite{MikhailetsMurach14,AnopDenkMurach21CPAA,MikhailetsMurach12BJMA2}.

В роботах \cite{DiachenkoLos22JEPE,DiachenkoLos23UMJ8} доведено теореми про ізоморфізми, досліджено регулярність та класичність узагальнених розв'язків
загальної параболічної крайової задачі для системи диференціальних рівнянь другого порядку в шкалі узагальнених просторів Соболєва.
При цьому в означенні класичного розв'язку не вимагалась його неперервність на лінії з'єднання основи і бічної поверхні циліндра.
Часто зазначену умову неперервності вважають частиною означення класичного розв'язку задачі (див., наприклад, \cite[с.42]{MathEncikl-5}).
В цьому випадку класичний розв'язок називатимемо сильно класичним. Представляє інтерес отримати умови, за яких узагальнений розв'язок параболічної задачі для систем рівнянь
буде сильно класичним. В скалярному випадку такий результат отримано в \cite{Los16UMJ9}.

Крайові задачі для параболічних систем диференціальних рівнянь другого порядку є математичними моделями багатьох прикладних задач.
Розглянемо в цій роботі окремий змістовний випадок параболічної задачі. А саме, мішану задачу для системи двох диференціальних рівнянь другого
порядку з двома крайовими умовами, одна з яких Діріхле, друга -- Неймана.
Такі задачі виникають, зокрема, в теорії тепломасообміну \cite[п.2.4]{Eidelman94}.
Встановимо нові достатні умови, за яких узагальнений розв'язок задачі буде сильно класичним, а також конкретизуємо умови існування класичного розв'язку цієї задачі
на основі результатів роботи \cite{DiachenkoLos23UMJ8}. Результати сформулюємо в термінах належності правих частин задачі узагальненим анізотропним просторам Соболєва.

\textbf{2. Постановка задачі.}
Нехай довільно задані ціле
число $n\geq2$, дійсне число $\tau>0$
і обмежена область $G\subset\mathbb{R}^{n}$ з
нескінченно гладкою межею $\Gamma:=\partial G$.
Позначимо $\Omega:=G\times(0,\tau)$ --- відкритий циліндр в $\mathbb{R}^{n+1}$,
$S:=\Gamma\times(0,\tau)$~--- його бічна поверхня.
Тоді $\overline{\Omega}:=\overline{G}\times[0,\tau]$ і
$\overline{S}:=\Gamma\times[0,\tau]$ є замикання
$\Omega$ і $S$ відповідно. Будемо ототожнювати $G$ з нижньою основою $G\times\{0\}$ циліндра $\Omega$.
Для частинних похідних функції, яка залежить від $x=(x_1,\ldots,x_n)\in\mathbb{R}^{n}$ і $t\in\mathbb{R}$ будемо використовувати наступні позначення:
$D^\alpha_x:=D^{\alpha_1}_{1}\dots D^{\alpha_n}_{n}$, де $D_{k}:=\partial/\partial{x_k}$ і $\partial_{t}:=\partial/\partial t$.
Тут $\alpha=(\alpha_1,...,\alpha_n)$ мультиіндекс, і $|\alpha|:=\alpha_1+\cdots+\alpha_n$.

Розглянемо  у циліндрі $\Omega$ таку систему диференціальних рівнянь:
\begin{equation}\label{28f1}
\begin{split}
\partial_{t} u_1(x,t)=
&a_{11}\Delta u_{1}(x,t)+a_{12}\Delta u_{2}(x,t)+f_{1}(x,t),\\
\partial_{t} u_2(x,t)=
&a_{21}\Delta u_{1}(x,t)+a_{22}\Delta u_{2}(x,t)+f_{2}(x,t),\\
&\mbox{для всіх}\quad (x,t)\in\Omega.
\end{split}
\end{equation}
Тут всі коефіцієнти $a_{ij}$ є сталі дійсні числа, а характеристичні числа $\lambda_1$ і $\lambda_2$ матриці $(a_{ij})$ такі, що $\lambda_2>\lambda_1>0$.
Така система є параболічною за Петровським (див., наприклад, \cite[п.2.4]{Eidelman94}).

На бічній поверхні циліндра задано дві крайові умови:
\begin{equation}\label{28f2}
\begin{split}
b_{11}\,D_{n}&u_{1}(x,t)+b_{12}\,D_{n} u_{2}(x,t)=g_{1}(x,t),\\
b_{21}&u_{1}(x,t)+b_{22} u_{2}(x,t)=g_{2}(x,t),\\
&\mbox{для всіх}\quad (x,t)\in S.
\end{split}
\end{equation}
Тут всі коефіцієнти $b_{ij}$ є сталі дійсні числа такі, що для крайових операторів виконуються так звані умови доповнюваності \cite[п.2.4]{Eidelman94}, необхідні для коректної постановки задачі.
Ці умови виконуються, зокрема, якщо $b_{11}b_{22}=0$ або $b_{21}b_{12}=0$, але не одночасно.

На основі циліндра задано початкові дані Коші:
\begin{equation}\label{28f3}
\begin{split}
&u_{1}(x,t)\big|_{t=0}=h_{1}(x),\\
&u_{2}(x,t)\big|_{t=0}=h_{2}(x),\\
&\mbox{для всіх}\quad x\in G.
\end{split}
\end{equation}

Початково-крайова задача
\eqref{28f1}--\eqref{28f3} буде параболічною за Петровським
у циліндрі $\Omega$ (див. означення в \cite[розд.~1, \S~1]{Solonnikov65} або \cite{Eidelman94}).
Дійсно, як зазначили вище, система \eqref{28f1} є параболічною за Петровським, а коефіцієнти крайових операторів \eqref{28f2} такі,
що для останніх виконується умова доповнюваності.

\textbf{3. Функціональні простори.}
Для зручності читання роботи  нагадаємо коротко означення узагальнених анізотропних просторів Соболєва, які потрібні для формулювання результатів.
Для викладу цього пункту скористаємось \cite[п.2]{DiachenkoLos23UMJ8} та \cite[п.4]{LosMikhailetsMurach21CPAA}.

Через $\mathcal{M}$ позначимо клас усіх неперервних функцій
$\varphi:[1,\infty)\rightarrow(0,\infty)$ таких, що:

(i) $\varphi$ і $1/\varphi$ обмежені на кожному відрізку $[1,c]$, де $1<c<\infty$;

(ii) $\varphi$ повільно змінна за Й.~Карамата на нескінченності, тобто
\begin{equation*}
\lim_{r\rightarrow\infty}\frac{\varphi(\lambda r)}{\varphi(r)}=1\quad\mbox{для кожного}\quad\lambda>0.
\end{equation*}

Нехай $s\in\mathbb{R}$ і $\varphi\in\mathcal{M}$.
За означенням, комплексний лінійний простір $H^{s,s/2;\varphi}(\mathbb{R}^{k})$, де $2\leq k\in\mathbb{Z}$, складається з усіх повільно зростаючих розподілів $w\in \mathcal{S}'(\mathbb{R}^{k})$ таких, що їх (повне) перетворення
Фур'є $\widetilde{w}$ є функцією, яка локально інтегровна на $\mathbb{R}^{k}$ за Лебегом і задовольняє умову
\begin{equation}\label{norm}
\begin{split}
&\|w\|_{H^{s,s/2;\varphi}(\mathbb{R}^{k})}:=\\
&\biggl(\;\int\limits_{\mathbb{R}^{k-1}}\int\limits_{\mathbb{R}}
\bigl(1+|\xi|^2+|\eta|\bigr)^{s}
\,\varphi^{2}\bigl((1+|\xi|^2+|\eta|)^{1/2}\bigr)\,
|\widetilde{w}(\xi,\eta)|^{2}\,d\xi\,d\eta\biggr)^{1/2}<\infty,
\end{split}
\end{equation}
де $\xi\in\mathbb{R}^{k-1}$ і $\eta\in\mathbb{R}$.
Цей простір гільбертовий і сепарабельний відносно норми \eqref{norm}.
Він є окремим випадком просторів $\mathcal{B}_{p,\mu}$, введених Л. Хермандером \cite[п.~2.2]{Hermander63}; а саме, $H^{s,s/2;\varphi}(\mathbb{R}^{k})=
\mathcal{B}_{p,\mu}$ за умови, що $p=2$ і
$$\mu(\xi,\eta)\equiv \bigl(1+|\xi|^2+|\eta|\bigr)^{s/2}\,\varphi\bigl((1+|\xi|^2+|\eta|)^{1/2}\bigr).
$$

Гільбертовий анізотропний простір $H^{s,s/2;\varphi}(\Omega)$ означається як простір звужень на $\Omega$ усіх розподілів з $H^{s,s/2;\varphi}(\mathbb{R}^{n+1})$, а
гільбертовий анізотропний простір $H^{s,s/2;\varphi}(S)$ на бічній поверхні циліндра означається за базовим простором $H^{s,s/2;\varphi}(\mathbb{R}^{n})$ за допомогою спеціальних
локальних карт на бічній поверхні циліндра (див. \cite[п.~1]{Los16JMathSci}).
Означення та основні властивості просторів $H^{s,s/2;\varphi}(W)$, де $W\in\{\Omega, S\}$, наведені, наприклад, в \cite[п.2]{DiachenkoLos22JEPE}.
Ізотропні простори $H^{s;\varphi}(G)$ і $H^{s;\varphi}(\Gamma)$, задані на основі $G$ циліндра та лінії $\Gamma$ з'єднання основи і бічної поверхні відповідно, означено в
\cite[п.2.1.1, 3.2.1]{MikhailetsMurach14}, \cite{MikhailetsMurach12BJMA2}.

Також буде потрібний простір правих частин задачі $\mathcal{Q}^{s-2,s/2-1;\varphi}$ (див. \cite[п.4]{DiachenkoLos22JEPE}).
Він утворений такими векторами $(f_1,f_2,g_1,g_2,h_1,h_2)$ з простору
\begin{align*}
\bigl(H^{s-2,s/2-1;\varphi}(\Omega)\bigr)^2&\oplus H^{s-3/2,s/2-3/4;\varphi}(S)\\
&\oplus H^{s-1/2,s/2-1/4;\varphi}(S)\oplus\bigl(H^{s-1;\varphi}(G)\bigr)^2,
\end{align*}
що задовольняють природні умови узгодження правих частин параболічної задачі \eqref{28f1}--\eqref{28f3}.

Нарешті, нагадаємо потрібні означення локальних версій просторів, про які йшла мова вище (див. \cite[п.4]{LosMikhailetsMurach21CPAA}).
Нехай $U$~--- відкрита множина в $\mathbb{R}^{n+1}$ така, що $\nobreak{\Omega_0:=U\cap\Omega\neq\varnothing}$ і $U\cap\Gamma=\varnothing$.
Покладемо $\Omega':=U\cap\partial\overline{\Omega}$, $S_0:=U\cap S$, $S':=U\cap \{(x,\tau):x\in\Gamma\}$ і $G_0:=U\cap G$.
Позначимо через $H^{s,s/2;\varphi}_{\mathrm{loc}}(\Omega_0,\Omega')$ лінійний простір усіх розподілів $u$ на $\Omega$ таких, що $\chi u\in H^{s,s/2;\varphi}(\Omega)$ для кожної функції $\chi\in C^\infty (\overline\Omega)$, яка задовольняє умову $\mathrm{supp}\,\chi\subset\Omega_0\cup\Omega'$.
Аналогічно, позначимо через $H^{s,s/2;\varphi}_{\mathrm{loc}}(S_0,S')$ лінійний простір усіх розподілів $v$ на $S$ таких, що $\chi v\in H^{s,s/2;\varphi}(S)$ для будь-якої функції $\chi\in C^\infty (\overline S)$, яка задовольняє умову $\mathrm{supp}\,\chi\subset S_0\cup S'$. Нарешті, $H^{s;\varphi}_{\mathrm{loc}}( G_0)$ позначає лінійний простір усіх розподілів $w$ на $G$ таких, що $\chi w\in H^{s;\varphi}(G)$ для кожної функції $\chi\in C^\infty (\overline G)$, яка задовольняє умову $\mathrm{supp}\,\chi\subset G_0$.

Якщо $\varphi(\cdot)=1$, то ці простори є соболєвськими.
В цьому випадку прибираємо індекс $\varphi$ у їх позначеннях.

\textbf{4. Основні результати.}
Параболічність за Петровським задачі \eqref{28f1}--\eqref{28f3} означає її коректну розв'язність у відповідних шкалах узагальнених просторів Соболєва
(див. \cite[Теорема~4.1]{DiachenkoLos22JEPE}).
Це, зокрема, означає, що для будь-якого вектора $(f_1,f_2,g_1,g_2,h_1,h_2)$ правих частин задачі з соболєвського простору
$\mathcal{Q}^{0,0}$ задача має єдиний розв'язок $(u_1,u_2)\in \bigl(H^{2,1}(\Omega)\bigr)^2$.

В нашому випадку простір $\mathcal{Q}^{0,0}$ складається з вектор-функцій $(f_1,f_2,g_1,g_2,h_1,h_2)$, які належать простору
\begin{equation*}
\bigl(H^{0,0}(\Omega)\bigr)^2\oplus H^{1/2,1/4}(S)\oplus H^{3/2,3/4}(S)
\oplus\bigl(H^{1}(G)\bigr)^2
\end{equation*}
і задовольняють природню умову  узгодження на лінії $\Gamma$ з'єднання бічної поверхні і основи циліндра:
\begin{equation}\label{28f-copmcond}
g_2\!\upharpoonright\!\Gamma=(b_{21}h_{1}+b_{22}h_{2})\!\upharpoonright\!\Gamma.
\end{equation}

Отже, для будь-якого вектора $(f_1,f_2,g_1,g_2,h_1,h_2)$ правих частин задачі \eqref{28f1}--\eqref{28f3}, що задовольняє такі умови
\begin{equation*}
\begin{split}
&f_1,f_2\in H^{0,0}(\Omega),\\
&g_1\in H^{1/2,1/4}(S),\,\,g_2\in H^{3/2,3/4}(S),\\
&h_1,h_2\in H^{1}(G),\\
&g_2\!\upharpoonright\!\Gamma=(b_{21}h_{1}+b_{22}h_{2})\!\upharpoonright\!\Gamma,
\end{split}
\end{equation*}
задача має єдиний розв'язок $u=(u_1,u_2)\in \bigl(H^{2,1}(\Omega)\bigr)^2$.
Цей розв'язок називаємо \it узагальненим \rm розв'язком нашої задачі.

З практичної точки зору є важливим питання, за яких умов на праві частини задачі її узагальнений розв'язок буде в певному розумінні класичним.
Оскільки під класичним розуміють такий неперервно диференційовний певну кількість разів розв'язок $u=(u_1,u_2)$, що  ліві частини системи, крайових та початкових умов
обчислюються в сенсі класичного диференціювання та слідів неперервних функцій $u_1$ та $u_2$.
Як зазначали у вступі, часто однією з умов в означенні класичного розв'язку задачі є умова його неперервності на лінії
$\Gamma$ з'єднання бічної поверхні і основи циліндра. Ця умова буде виконана, якщо розв'язок є неперервним у замкненому циліндрі
$\overline\Omega$. В такому випадку будемо вживати термін "сильно класичний розв'язок".

Дамо точні означення.
Нехай
\begin{equation*}
S_{\varepsilon}:=\{x\in\Omega:\mbox{dist}(x,S)<\varepsilon\},\quad
G_{\varepsilon}:=\{x\in\Omega:\mbox{dist}(x,G)<\varepsilon\},
\end{equation*}
де число $\varepsilon>0$.

Узагальнений розв'язок $u=(u_1,u_2)\in \bigl(H^{2,1}(\Omega)\bigr)^2$ задачі \eqref{28f1}--\eqref{28f3}
назвемо \it сильно класичним, \rm якщо він та його узагальнені похідні задовольняють такі умови:
\begin{itemize}
\item [(a1)] $D^\alpha_xu_1$, $D^\alpha_xu_2$ при $0\leq|\alpha|\leq 2$ та $\partial_tu_1$ і $\partial_tu_2$ неперервні на $\Omega$;
\item [(b1)] $D^\alpha_x u_1$ і $D^\alpha_x u_2$ неперервні на $S_{\varepsilon}\cup S$ для деякого числа $\varepsilon>0$, якщо $0\leq|\alpha|\leq 1$;
\item [(c1)] $u_1$ і $u_2$ неперервні на $\overline\Omega$.
\end{itemize}

Нагадаємо, узагальнений розв'язок $u=(u_1,u_2)\in \bigl(H^{2,1}(\Omega)\bigr)^2$ задачі \eqref{28f1}--\eqref{28f3}
називається \it класичним, \rm (див. \cite[п.3]{DiachenkoLos23UMJ8}) якщо він та його узагальнені похідні задовольняють умови (a1), (b1) та
\begin{itemize}
\item [(c2)] $u_1$ і $u_2$ неперервні на $G_{\varepsilon}\cup G$ для деякого числа $\varepsilon>0$.
\end{itemize}

В означенні класичного розв'язку задачі формулюються мінімальні умови на вектор-функцію $u$, за яких вона в термінах класичних похідних  і слідів
задовольняє рівняння \eqref{28f1}, крайові умови \eqref{28f2} і початкові умови \eqref{28f3}. Для цього умова (c1) неперервності
вектор-функції $u$ в усьому циліндрі  заміняється на умову (c2) неперервності лише в деякому малому околі основи $G$.

Якщо вектор-функція $u=u(x,t)=(u_1(x,t), u_2(x,t))$ є класичним розв'язком задачі \eqref{28f1}--\eqref{28f3}, то  ліві частини рівняння, крайових та початкових умов є неперервними функціями на відповідних множинах. Зрозуміло, що сильно класичний розв'язок є класичним, але не навпаки.

Умови, за яких узагальнений розв'язок загальної параболічної крайової задачі для системи диференціальних рівнянь є класичним сформульовані в \cite[Теорема 4]{DiachenkoLos23UMJ8}.
Наступна теорема є конкретизацією зазначеної теореми для розглядуваної задачі.

\textbf{Теорема 1.} \it
Припустимо, що вектор-функція $u\in \bigl(H^{2,1}(\Omega)\bigr)^2$ є узагальненим розв'язком параболічної задачі \eqref{28f1}--\eqref{28f3}, праві частини якої задовольняють такі умови:
\begin{equation}\label{28f12}
\begin{aligned}
(f_1,f_2)\in\, &\bigl(H_{\mathrm{loc}}^{1+n/2,\,1/2+n/4;\varphi}
(\Omega,\varnothing)\bigr)^2\cap\\
&\cap
\bigl(H_{\mathrm{loc}}^{n/2,\,n/4;\varphi}
(S_{\varepsilon},S)\bigr)^2\cap\\
&\cap
\bigl(H_{\mathrm{loc}}^{-1+n/2,\,-1/2+n/4;\varphi}
(G_{\varepsilon},G)\bigr)^2,\\
g_1\in H_{\mathrm{loc}}^{n/2+1/2,\,n/4+1/4;\varphi}
&(S,\varnothing),\,\,g_2\in H_{\mathrm{loc}}^{n/2+3/2,\,n/4+3/4;\varphi}
(S,\varnothing),\\
(h_1,h_2)&\in\bigl(H_{\mathrm{loc}}^{n/2;\varphi}(G)\bigr)^2
\end{aligned}
\end{equation}
з деяким функціональним параметром $\varphi\in\mathcal{M}$, що задовольняє інтегральну умову
\begin{equation}\label{IntCond}
\int\limits_{1}^{\,\infty}\;\frac{dr}{r\,\varphi^2(r)}<\infty.
\end{equation}
У випадку $n=2$ додатково припускаємо, що функція $\varphi$ зростає (в нестрогому сенсі).
Тоді розв'язок $u=(u_1,u_2)$ класичний.
\rm

Теорему~1 отримуємо з \cite[Теорема 4]{DiachenkoLos23UMJ8} при $N=2$, $l_0=l_1=1$ і $l_2=0$.

\medskip

Відмітимо, що припущення зростання $\varphi$ при $n=2$ необхідне нам з огляду на те, що результат про розв'язність задачі
\eqref{28f1}--\eqref{28f3} в узагальнених просторах Соболєва \cite[Теорема~4.1]{DiachenkoLos22JEPE} отримано за умови, що її праві частини належать просторам, які вкладаються у відповідні простори $L_2(\cdot)$.
А, наприклад, простір $H^{0,0;\varphi}(\Omega)$ для спадної $\varphi$ буде вже ширший ніж $H^{0,0}(\Omega)=L_2(\Omega)$.

{\bf Зауваження 1}. Використання узагальнених просторів Соболєва дає більш тонкий результат ніж у випадку просторів Соболєва.
А саме, щоб висновок теореми 1 залишився правильним при $\varphi=1$ (тут \eqref{IntCond} не виконується), треба в \eqref{28f12} збільшити числові показники регулярності просторів
на деяке число $\delta>0$. Іншими словами ці умови замінити на такі:
\begin{equation*}\label{28f12-s}
\begin{aligned}
(f_1,f_2)\in\, &\bigl(H_{\mathrm{loc}}^{1+n/2+\delta,\,1/2+n/4+\delta/2}
(\Omega,\varnothing)\bigr)^2\cap\\
&\cap
\bigl(H_{\mathrm{loc}}^{n/2+\delta,\,n/4+\delta/2}
(S_{\varepsilon},S)\bigr)^2\cap\\
&\cap
\bigl(H_{\mathrm{loc}}^{-1+n/2+\delta,\,-1/2+n/4+\delta/2}
(G_{\varepsilon},G)\bigr)^2,
\end{aligned}
\end{equation*}
\begin{equation*}\label{28f12bis-s}
g_1\in H_{\mathrm{loc}}^{n/2+1/2+\delta,\,n/4+1/4+\delta/2}
(S,\varnothing),\,\,g_2\in H_{\mathrm{loc}}^{n/2+3/2+\delta,\,n/4+3/4+\delta/2}
(S,\varnothing),
\end{equation*}
\begin{equation*}\label{28f9-s}
(h_1,h_2)\in\bigl(H_{\mathrm{loc}}^{n/2+\delta}(G)\bigr)^2.
\end{equation*}
Останні умови є більш сильними ніж \eqref{28f12}, оскільки для довільних $s_1>s_2$ і $\varphi\in\mathcal{M}$
правильні вкладання $H^{s_1,s_1/2}(\cdot)\hookrightarrow H^{s_2,s_2/2;\varphi}(\cdot)$ \cite[c. 22]{LosMikhailetsMurach21Monograph}.

Сформулюємо тепер основний результат -- умови, за яких розв'язок $u=u(x,t)$ задачі \eqref{28f1}--\eqref{28f3} буде сильно класичним.

\medskip

\textbf{Теорема 2.} \it
Припустимо, що вектор-функція $u\in \bigl(H^{2,1}(\Omega)\bigr)^2$ є узагальненим розв'язком параболічної задачі \eqref{28f1}--\eqref{28f3}, праві частини якої задовольняють такі умови:
\begin{equation}\label{28f12-sk}
(f_1,f_2)\in\, \bigl(H_{\mathrm{loc}}^{1+n/2,\,1/2+n/4;\varphi}
(\Omega,\varnothing)\bigr)^2\cap
\bigl(H_{\mathrm{loc}}^{n/2,\,n/4;\varphi}
(S_{\varepsilon},S)\bigr)^2
\end{equation}
\begin{equation}\label{28f12bis-sk}
g_1\in H_{\mathrm{loc}}^{n/2+1/2,\,n/4+1/4;\varphi}(S,\varnothing),\quad
g_2\in H_{\mathrm{loc}}^{n/2+3/2,\,n/4+3/4;\varphi}(S,\varnothing),
\end{equation}
\begin{equation}\label{28f9-sk}
(f_1,f_2,g_1,g_2,h_1,h_2)\in\mathcal{Q}^{-1+n/2,-1/2+n/4;\varphi}
\end{equation}
з деяким функціональним параметром $\varphi\in\mathcal{M}$, що задовольняє інтегральну умову
\eqref{IntCond}.
У випадку $n=2$ додатково припускаємо, що функція $\varphi$ зростає (в нестрогому сенсі).
Тоді розв'язок $u=(u_1,u_2)$ сильно класичний.
\rm

Відмітимо, щоб вказати, яким просторам належать праві частини задачі згідно умови \eqref{28f9-sk}, треба скористатись означенням простору
$\mathcal{Q}^{-1+n/2,-1/2+n/4;\varphi}$. Він має доволі складну будову \cite[п.4]{DiachenkoLos22JEPE} завдяки умовам узгодження, накладеним на компоненти його елементів. Кількість та складність
цих умов зростає зі зростанням $n$.

{\bf Зауваження 2}. У практично важливих випадках $n=2$ або $n=3$ умова \eqref{28f9-sk} в теоремі~2 еквівалентна такій
\begin{equation}\label{28f15}
\begin{aligned}
&(f_1,f_2)\in \bigl(H^{-1+n/2,\,-1/2+n/4;\varphi}(\Omega )\bigr)^2,\\
&g_1\in H^{n/2-1/2,\,n/4-1/4;\varphi}(S), \\
&g_2\in H^{n/2+1/2,\,n/4+1/4;\varphi}(S),\\
&(h_1,h_2)\in\bigl(H^{n/2;\varphi}(G)\bigr)^2.
\end{aligned}
\end{equation}
Взагалі, простір $\mathcal{Q}^{-1+n/2,-1/2+n/4;\varphi}$ для $n=2$ або $n=3$ складається з вектор-функцій, що задовольняють одночасно
\eqref{28f15} і \eqref{28f-copmcond}.
Для вектор-функцій $(f_1,f_2,g_1,g_2,h_1,h_2)$ з формулювання теореми~2 виконання \eqref{28f-copmcond} випливає з умови теореми $u\in \bigl(H^{2,1}(\Omega)\bigr)^2$. Дійсно, тоді $(f_1,f_2,g_1,g_2,h_1,h_2)\in\mathcal{Q}^{0,0}$. За означенням, елементи простору $\mathcal{Q}^{0,0}$ задовольняють \eqref{28f-copmcond}.
При $n\geq4$ простір $\mathcal{Q}^{-1+n/2,-1/2+n/4;\varphi}$ складається з векторів, які задовольняють включення \eqref{28f15} та більш
складні ніж \eqref{28f-copmcond} умови узгодження \cite[п.4]{DiachenkoLos22JEPE}.

Зазначимо, що для теореми~2 правильна версія зауваження~1.

{\bf Доведення теореми~2.} Покажемо, що функції $u_1$ і $u_2$ задовольняють умови (a1)--(c1) означення сильно класичного розв'язку.
Для умов (а1) і (b1) використаємо теорему~3 з \cite{DiachenkoLos23UMJ8}.

Почнемо з (а1). Покладемо $N=2$, $\Omega_0=\Omega$, $\Omega'=S_0=S'=G_0=\varnothing$ і $p=2$. Далі скористаємось включенням
$$
(f_1,f_2)\in\,\bigl(H_{\mathrm{loc}}^{1+n/2,\,1/2+n/4;\varphi}
(\Omega,\varnothing)\bigr)^2
$$
з умови~\eqref{28f12-sk}. Тоді з теореми~3 \cite{DiachenkoLos23UMJ8} випливає, що вектор-функція
$u(x,t)=(u_1(x,t),u_2(x,t))$ і кожна її узагальнена частинна похідна
$D_{x}^{\alpha}\partial_{t}^{\beta}u(x,t)=(D_{x}^{\alpha}\partial_{t}^{\beta}u_1(x,t),D_{x}^{\alpha}\partial_{t}^{\beta}u_2(x,t))$,
де $|\alpha|+2\beta\leq 2$,  неперервні на множині $\Omega_0\cup\Omega'=\Omega$.
При $\beta=0$ маємо неперервність $D^\alpha_xu_1$ і $D^\alpha_xu_2$ для $|\alpha|\leq 2$, при  $|\alpha|=0$ -- неперервність $\partial_tu_1$ і $\partial_tu_2$ відповідно.
Таким чином, виконується умова (а1).

Перейдемо до умови (b1). Покладемо $N=2$, $\Omega_0=S_\varepsilon$, $\Omega'=S_0=S$, $S'=G_0=\varnothing$ і $p=1$.
Використаємо включення
$$
(f_1,f_2)\in \bigl(H_{\mathrm{loc}}^{n/2,\,n/4;\varphi}
(S_{\varepsilon},S)\bigr)^2
$$
з \eqref{28f12-sk} та включення \eqref{28f12bis-sk}:
$$
g_1\in H_{\mathrm{loc}}^{n/2+1/2,\,n/4+1/4;\varphi}(S,\varnothing),\quad
g_2\in H_{\mathrm{loc}}^{n/2+3/2,\,n/4+3/4;\varphi}(S,\varnothing).
$$
В цьому випадку з теореми~3 \cite{DiachenkoLos23UMJ8} випливає, що вектор-функція
$u(x,t)=(u_1(x,t),u_2(x,t))$ і кожна її узагальнена частинна похідна
$D_{x}^{\alpha}\partial_{t}^{\beta}u(x,t)=(D_{x}^{\alpha}\partial_{t}^{\beta}u_1(x,t),D_{x}^{\alpha}\partial_{t}^{\beta}u_2(x,t))$,
де $|\alpha|+2\beta\leq 1$,  неперервні на множині $\Omega_0\cup\Omega'=S_\varepsilon\cup S$.
При $\beta=0$ маємо неперервність $D^\alpha_xu_1$ і $D^\alpha_xu_2$ для $|\alpha|\leq 1$.
Отже, виконується умова (b1).

Нарешті, розглянемо умову (c1).
Скориставшись включенням \eqref{28f9-sk}
$$(f_1,f_2,g_1,g_2,h_1,h_2)\in\mathcal{Q}^{-1+n/2,-1/2+n/4;\varphi}$$
та \cite[Теорема 1]{DiachenkoLos23UMJ8} з $s=1+n/2$ робимо висновок, що
$$(u_1(x,t),u_2(x,t))\in (H^{1+n/2,1/2+n/4;\varphi}(\Omega))^2.$$
Розглянемо простір $H^{1+n/2,1/2+n/4;\varphi}(\mathbb{R}^{n+1})$ з функціональним параметром $\varphi$, який задовольняє інтегральну умову \eqref{IntCond}.
З \cite[Теорема~1.13(і)]{LosMikhailetsMurach21Monograph} при $p=0$ і $b=1$ випливає, що всі елементи цього простору є неперевними на $\mathbb{R}^{n+1}$ функціями.
За означенням, простір $H^{1+n/2,1/2+n/4;\varphi}(\Omega)$ складається зі звужень на $\Omega$ розподілів з простору $H^{1+n/2,1/2+n/4;\varphi}(\mathbb{R}^{n+1})$.
Тому для $u_1$ і $u_2$ існують такі $w_1$ і $w_2$ з простору $H^{1+n/2,1/2+n/4;\varphi}(\mathbb{R}^{n+1})$, що $u_1=w_1\!\upharpoonright\!\Omega$ і $u_2=w_2\!\upharpoonright\!\Omega$. Це означає, що $u_1$ і $u_2$ є неперервними функціями на $\overline\Omega$, тобто виконується умова (c1).

Теорему 2 доведено.

\end{document}